\documentclass[a4paper,12pt]{article}
\usepackage{graphicx}
\usepackage{epsfig}
\usepackage{fancyhdr}

\usepackage{amsmath,amsfonts,amssymb,amsthm}
\usepackage{setspace}
\newtheorem{theorem}{Theorem}[section]

\newcommand{\reftit}{\textit}    
\newcommand{\refis}{\textbf}     

\begin{document}

\title{Occupation times via Bessel functions.}
\author{Yevgeniy Kovchegov\footnote{Department of Mathematics,
 Oregon State University, Corvallis, OR  97331-4605, USA},
  Nick Meredith\footnote{University of California, Berkeley, CA 94720, USA} and Eyal Nir\footnote{Department of Chemistry and Biochemistry, University of California, Los Angeles, CA 90095, USA} \footnote{E.N. is supported by the Human Frontier Science Program (HFSP)}\\
}
\date{ }
\maketitle

{\abstract 
This study of occupation time densities for continuous-time Markov processes was inspired by the work of E.Nir et al (see \cite{nir}) in the field of Single Molecule FRET spectroscopy. There, a single molecule fluctuates between two or more states, and the experimental observable depends on the state's occupation time distribution. To mathematically describe the observable there was a need to calculate a single state occupation time distribution.\\
\indent In this paper, we consider a Markov process with countably many states. In order to find a one-stete occupation time density,
we use a combination of Fourier and Laplace transforms in the way that allows for inversion of the Fourier transform.
We derive an explicit expression for an occupation time density in the case of a simple continuous time random walk on $\mathbb{Z}$.
Also we examine the spectral measures in Karlin-McGregor diagonalization in an attempt to represent occupation 
time densities via modified Bessel functions.} 
\vskip 0.3 in
\noindent{\bf Primary Subjects:} 60J35, 47N30, 47N20, 33C90, 33C10



\section{Introduction}

 The occupation time densities for continuous-time Markov processes that live on countable state space was a subject
  of intense research in the 60s, 70s and early 80s.  We would like to refer the reader to   \cite{d}, \cite{dm}, \cite{kovaleva},
  \cite{p}, \cite{w} and \cite{b} (many published in the Journal of Applied Probability) for some of the results in the field.
  With the exception of \cite{b}, the main instrument was the multidimensional Laplace transform. Also, the reader can
   check \cite{s} and \cite{bmns} for the most recent developments in the field.

 This paper is a mathematical follow-up to the research done by E.Nir et al \cite{nir} in the field of Single Molecule Fluorescence Resonance Energy Transfer (FRET) spectroscopy, where a single molecule fluctuates between two or more states, and the experimental observable depends on the state's occupation time distribution. While working on \cite{nir} the authors have noticed that the single state occupation time densities, when computed via randomization
 technique (i.e. multiple infinite sums) can often be represented via modified Bessel functions of the kind
 $$I_{n}(z)=\sum_{k=0}^{\infty}{\frac{1}{k!(k+n)!}\left(\frac{z}{2}\right)^{2k+n}}~.$$
 In this paper, we use spectral theory in an attempt to find a rigorous explanation for the relationship between occupation times and Bessel functions.  We will show a connection between a spectral measure of a generator and a Laplace transform of 
 a single state occupation time distribution taken with respect to a time variable $t$.

  The occupation times for birth-and-death chains were studied by Karlin and McGregor with
orthogonal polynomials in \cite{km3} following the paper of Darling and Kac, see \cite{dk} .
 Both papers considered occupation times for Markov processes when $t$ is taken to $\infty$, while this paper concerns
 with explicit expressions for a given fixed time interval $[0,t]$.

\section{Approach and results}

 In this section we will state and prove theorems and formulas that will relate occupation times to the spectral measure of a generator,
 we will find an explicit expression for a one-dimensional symmetric random walk (see Theorem 2.2), and find a representation
 of occupation time distributions via modified Bessel functions.

\subsection{General case: spectral representation}

 Consider an irreducible discrete Markov process with the generator matrix (or operator, if the state space is infinite).
 $$
Q =
\begin{pmatrix}
 -\sum_{j \neq 0} \lambda_{0,j}& \lambda_{0,1} & ... 
 \\
 \\
 \lambda_{1,0}& -\sum_{j \neq 1} \lambda_{1,j} & ... 
 \\
 \\
 ...& ... & ... 
\end{pmatrix}
$$
 Let $\{0,1,\dots \}$ be a countable state space, and for a time interval $[0,t]$, let $f_k(t,x)$ denote 
 the probability density function for the occupation time associated with state $0$, given that the continuous-time process 
 commences at state $k$.
 Denote $e_0=\left[\begin{array}{c}1\\ 0\\ 0\\ \vdots \end{array}\right]$.
 
 \noindent
 First, we derive the following expression for the occupation time density $f_0(t,x)$.
 
 {\theorem
  The Laplace transform w.r.t. time t of $f_0(t,x)$ can be written as
  \begin{eqnarray} \label{generalLaplace}
   {\L}_{f_0}(s_1,x)={1 \over s_1 h(s_1)}\exp\left\{-{x  \over h(s_1)}\right\}~,
 \end{eqnarray}
 where $h(s)=-\left((Q-sI)^{-1}e_0, e_0\right)$.
 }
 
 \begin{proof}
 The integral equations relating $\{f_k(t,x)\}_{k=0,1,\dots}$ can be produced via conditioning as follows:
 \begin{eqnarray*}
 f_{0}(t,x) & = &
e^{-(\sum_{m:~m\not=1} \lambda_{0,m})t}\delta_{t}(x) +
\sum_{k:~k\not=0}\int_{0}^{t}f_{k}(t-y,x-y){\lambda_{0,k}}e^{-(\sum_{m:~m\not=0} \lambda_{0,m})y}dy ,\\
 f_{j}(t,x) & = &
e^{-(\sum_{m:~m\not=j} \lambda_{j,m})t}\delta_{0}(x) +
\sum_{k:~k\not=j}\int_{0}^{t}f_{k}(t-y,x){\lambda_{j,k}}e^{-(\sum_{m:~m\not=j} \lambda_{j,m})y}dy \\
&  & \text{ for }j=1,2,\dots ~. 
 \end{eqnarray*}
 Plugging in $\psi = t-y$ into the above equation, we obtain
  \begin{eqnarray*}
 f_{0}(t,x) & = &
e^{-(\sum_{m:~m\not=0} \lambda_{0,m})t}\delta_{t}(x) +
\sum_{k:~k\not=0}\int_{0}^{t}f_{k}(\psi,x-t+\psi){\lambda_{0,k}}e^{-(\sum_{m:~m\not=0} \lambda_{0,m})(t-\psi)}d\psi ,\\
 f_{j}(t,x) & = &
e^{-(\sum_{m:~m\not=j} \lambda_{j,m})t}\delta_{0}(x) +
\sum_{k:~k\not=j}\int_{0}^{t}f_{k}(\psi,x){\lambda_{j,k}}e^{-(\sum_{m:~m\not=j} \lambda_{j,m})(t-\psi)}d\psi \\
&  & \text{ for }j=1,2,\dots ~. 
 \end{eqnarray*}
 Now taking the Fourier transform with respect to $x$ we arrive at
  \begin{eqnarray*}
 \hat{f}_{0}(t,s_2) \!\!\!\! & = & \!\!\!\! e^{-(\sum_{m:~m\not=0} \lambda_{0,m} -is_2)t}+
\sum_{k:~k\not=0}\int_{-\infty}^{\infty} \int_{0}^{t}f_{k}(\psi,x-t+\psi){\lambda_{0,k}}e^{-(\sum_{m:~m\not=0} \lambda_{0,m})(t-\psi)}d\psi e^{is_2x}dx,\\
 \hat{f}_{j}(t,s_2) \!\!\!\! & = & \!\!\!\! e^{-(\sum_{m:~m\not=j} \lambda_{j,m})t} +
\sum_{k:~k\not=j}\int_{-\infty}^{\infty} \int_{0}^{t}f_{k}(\psi,x){\lambda_{j,k}}e^{-(\sum_{m:~m\not=j} \lambda_{j,m})(t-\psi)}d\psi e^{is_2x}dx\\
&  & \text{ for }j=1,2,\dots ~.
 \end{eqnarray*}
  \noindent The above equations simplify to
  \begin{eqnarray*}
 e^{(\sum_{m:~m\not=0} \lambda_{0,m} -is_2)t}\hat{f}_{0}(t,s_2)  & = &  1+
\sum_{k:~k\not=0} \int_{0}^{t}\hat{f}_{k}(\psi,s_2){\lambda_{0,k}}e^{(\sum_{m:~m\not=0} \lambda_{0,m}-is_2)\psi}d\psi ,\\
 e^{(\sum_{m:~m\not=j} \lambda_{j,m})t} \hat{f}_{j}(t,s_2)  & = &  1 +
\sum_{k:~k\not=j} \int_{0}^{t}\hat{f}_{k}(\psi,s_2){\lambda_{j,k}}e^{(\sum_{m:~m\not=j} \lambda_{j,m})\psi}d\psi \\
&  & \text{ for }j=1,2,\dots ~.
 \end{eqnarray*}  
 Differentiating w.r.t. variable $t$, obtain
 \begin{eqnarray*}
 \left(\sum_{m:~m\not=0} \lambda_{0,m} -is_2 \right)\hat{f}_{0}(t,s_2)+{\partial \over \partial t} \hat{f}_{0}(t,s_2)  & = &  
\sum_{k:~k\not=0}{\lambda_{0,k}} \hat{f}_{k}(t,s_2),\\
 \left(\sum_{m:~m\not=j} \lambda_{j,m}\right)\hat{f}_{j}(t,s_2)+{\partial \over \partial t} \hat{f}_{j}(t,s_2)  & = &  
\sum_{k:~k\not=j}{\lambda_{j,k}} \hat{f}_{k}(t,s_2) ~~~(j=1,2,\dots )~.
 \end{eqnarray*}  
 Observe: $\hat{f}_j(0,s_2)=1$ for all $j$. Our next step is to take the Laplace transform w.r.t. variable $t$:
  \begin{eqnarray*} 
 \left(\sum_{m:~m\not=0} \lambda_{0,m}+s_1 -is_2 \right){\L}_{\hat{f}_{0}}(s_1,s_2)  & = &  
1+\sum_{k:~k\not=0}{\lambda_{0,k}} {\L}_{\hat{f}_{k}}(s_1,s_2),\\ 
 \left(\sum_{m:~m\not=j} \lambda_{j,m}+s_1\right){\L}_{\hat{f}_{j}}(s_1,s_2) & = &  
1+\sum_{k:~k\not=j}{\lambda_{j,k}} {\L}_{\hat{f}_{k}}(s_1,s_2) ~~~(j=1,2,\dots )~. 
 \end{eqnarray*}   
 \vskip 0.3 in
 \noindent The above system of equations can be rewritten via the spectral decomposition of the generator operator $Q$ as follows. 
 Let ${\L}_{\hat{f}}(s_1,s_2)=\left[\begin{array}{c}{\L}_{\hat{f}_{0}}(s_1,s_2) \\ {\L}_{\hat{f}_{1}}(s_1,s_2) \\ \vdots \end{array}\right]$
 and ${\bf 1}=\left[\begin{array}{c}1\\ 1\\ \vdots \end{array}\right]$. So we proved the following spectral identity 
 \begin{eqnarray} \label{generalFourierLaplace}
 (Q-s_1I){\L}_{\hat{f}}(s_1,s_2)=-\mathbf{1}-is_2\left[\begin{array}{c}{\L}_{\hat{f}_{0}}(s_1,s_2)\\ 0\\ 0\\ \vdots \end{array}\right]~.
 \end{eqnarray}
\vskip0.3 in
 \noindent Thus
 ${\L}_{\hat{f}}(s_1,s_2)=-(Q-s_1I)^{-1}\mathbf{1}-is_2{\L}_{\hat{f}_{0}}(s_1,s_2) (Q-s_1I)^{-1} e_0$ and
 $$ {\L}_{\hat{f}_{0}}(s_1,s_2)=-\left((Q-s_1I)^{-1}\mathbf{1}, e_0\right)-is_2 {\L}_{\hat{f}_{0}}(s_1,s_2) \left((Q-s_1I)^{-1} e_0,e_0 \right)~.$$
 Therefore the Laplace-Fourier transform of $f_0$ can be represented as
 \begin{eqnarray} \label{FourierLaplace1}
 {\L}_{\hat{f}_{0}}(s_1,s_2)={-\left((Q-s_1I)^{-1}\mathbf{1}, e_0\right) \over 1+is_2\left((Q-s_1I)^{-1} e_0,e_0 \right)}~.
 \end{eqnarray}
 Observe that ${\L}_{\hat{f}_{k}}(s_1,0)=\int_{[0,+\infty)} \int_{\mathbb{R}}e^{-s_1t} f_k(t,x) dx dt={1\over s_1}$ for all $k$.
 Substituting $s_2=0$ into (\ref{generalFourierLaplace}) gives
 $${1 \over s_1}(Q-s_1I)\mathbf{1}=-\mathbf{1}$$
 which is obviously true. This also implies
 $$(Q-s_1I)^{-1}\mathbf{1}=-{1 \over s_1}\mathbf{1}~.$$
 Therefore (\ref{FourierLaplace1}) can be simplified to
 \begin{eqnarray} \label{FourierLaplace}
 {\L}_{\hat{f}_{0}}(s_1,s_2)={1/s_1 \over 1-is_2h(s_1)}~,
 \end{eqnarray} 
 where $h(s)=-\left((Q-sI)^{-1}e_0, e_0\right)=\left((\int_0^{\infty} e^{-st} e^{Qt} dt)e_0,e_0\right)=\int_0^{\infty} e^{-st} p_{t}(0,0) dt$.
  The Fourier transform can be inverted via complex integration over a lower semi-circle contour with the radius
  converging to infinity:
  \begin{eqnarray} 
   {\L}_{f_0}(s_1,x)={1 \over s_1 h(s_1)}\exp\left\{-{x  \over h(s_1)}\right\}
   ={-1 \over s_1\!\left((Q-s_1I)^{-1} e_0,e_0 \right)} \exp\left\{{x  \over \left((Q-s_1I)^{-1} e_0,e_0 \right)}\right\}.
 \end{eqnarray}
 \end{proof}

\subsubsection{Example: Two-state Markov processes}

Consider a two-state Markov process with generator
$Q=\left(\begin{array}{cc}-\lambda & \lambda \\\mu & -\mu\end{array}\right)$.
Then
$$(Q-s_1I)^{-1}={-1 \over s_1^2+(\lambda+\mu)s_1} \left(\begin{array}{cc} \mu+s_1 & \lambda \\ \mu & \lambda+s_1\end{array}\right)$$
and (\ref{generalLaplace}) implies 
$${\L}(s_{1},x) 
 =  e^{-x(s_{1}+\lambda)}e^{\frac{{\lambda}{\mu}x}{s_{1}+\mu}}+
\frac{\lambda}{s_{1}+\mu}e^{-x(s_{1}+\lambda)}e^{\frac{{\lambda}{\mu}x}{s_{1}+\mu}}~.$$
Now, formula (29.3.81) of \cite{as} gives us  the following Laplace transforms
$$\int_{0}^{\infty}I_{0}(2\sqrt{at})e^{-pt}dt=\frac{1}{p}e^{\frac{a}{p}}~~
\text{ and }
~~\int_{0}^{\infty}\frac{1}{\sqrt{t}}I_{1}(2\sqrt{at})e^{-pt}dt=\frac{1}{\sqrt{a}}(e^{\frac{a}{p}}-1),$$
where $I_{0}(\cdot)$ and $I_{1}(\cdot)$ are modified Bessel functions.
Next, we rewrite the above identities as follows
$$e^{-px} \frac{1}{p}e^{\frac{a}{p}}=\int_{x}^{\infty}I_{0}(2\sqrt{a(t-x)})e^{-pt}dt$$
and
$$e^{-px} e^{\frac{a}{p}}=e^{-px}+\sqrt{a}\int_{x}^{\infty}\frac{1}{\sqrt{t-x}}I_{1}(2\sqrt{a(t-x)})e^{-pt}dt.$$
Let $a=\lambda \mu x$ and $p=s_1+\mu$. Plugging in, we get
$$\frac{1}{s_{1}+\mu}e^{-(s_{1}+\mu)x}e^{\frac{{\lambda}^{2}x}{s_{1}+\mu}}=e^{-px} \frac{1}{p}e^{\frac{a}{p}}
={\int}_{x}^{\infty}I_{0}(2\sqrt{\lambda \mu x(t-x)})e^{-\mu t}
e^{-s_1 t}dt$$ and thereforethe inverse Laplace transform of
$\frac{\lambda}{s_{1}+\mu}e^{-x(s_{1}+\lambda)}e^{\frac{{\lambda}{\mu}x}{s_{1}+\mu}}$
is $$\lambda e^{-\lambda x}e^{-\mu(t-x)}I_{0}(2\sqrt{{\lambda}{\mu}x(t-x)})$$ for $0 \leq x \leq t$.

Similarly
$$e^{-x(s_{1}+\mu)}e^{\frac{{\lambda}{\mu}x}{s_{1}+\mu}}=e^{-px} e^{\frac{a}{p}}
=e^{-\mu x} e^{-s_1 x}+\sqrt{{\lambda}{\mu}x}\int_{x}^{\infty}\frac{I_{1}(2\sqrt{{\lambda}{\mu}x(t-x)})}{\sqrt{t-x}}e^{-\mu t}e^{-s_{1}t}dt$$ 
which can be rewritten
as
$$e^{-x(s_{1}+\lambda)}e^{\frac{{\lambda}{\mu}x}{s_{1}+\mu}}=\sqrt{{\lambda}{\mu}x}\int_{x}^{\infty}\frac{I_{1}(2\sqrt{{\lambda}{\mu}x(t-x)})}{\sqrt{t-x}}e^{-{\lambda}x}e^{-\mu(t-x)}e^{-s_{1}t}dt+\int_{0}^{\infty}e^{-{\lambda}t}\delta_{t}(x)e^{-s_{1}t}dt$$
Therefore, the inverse Laplace transform of
$e^{-x(s_{1}+\lambda)}e^{\frac{{\lambda}{\mu}x}{s_{1}+\mu}}$ is
$$\sqrt{\frac{{\lambda}{\mu}x}{t-x}}I_{1}(2\sqrt{\lambda{\mu}x(t-x)})e^{-{\lambda}x}e^{-\mu(t-x)}+e^{-{\lambda}t}\delta_{t}(x)$$
 Here we do not divide by zero when $x=t$ as the $\sqrt{t-x}$ cancels on top and the bottom.
Adding the terms together, we obtain
$$f_{0}(t,x)=e^{-{\lambda}t}\delta_{t}(x)+\lambda
e^{-\lambda
x}e^{-\mu(t-x)}I_{0}(2\sqrt{{\lambda}{\mu}x(t-x)})+\sqrt{\frac{{\lambda}{\mu}x}{t-x}}I_{1}(2\sqrt{\lambda{\mu}x(t-x)})e^{-{\lambda}x}e^{-\mu(t-x)}$$
for $0 \leq x \leq t$.

The above equation was originally derived in \cite{p} via two-dimensional Laplace transform. One can also derive it via randomization,
 where the infinite sums are easily recognized to be the corresponding modified Bessel functions.

\subsubsection{Example: Three state Markov chain.}
 
  Here, the inversion of  ${\L}_{f_0}(s_1,x)$ as expressed in (\ref{generalLaplace}) with
  $${1 \over h(z)}=z-{\gamma_1 \over z+\beta_1}-{\gamma_2 \over z+\beta_2}$$
  can be expressed via convolutions of modified Bessel functions.
\vskip 0.2 in
\noindent The last two examples prompted us to look closely at the structure of the occupation time densities, and
 in an attempt to understand the mechanics of decomposing them via cylindrical functions $I_n(\cdot)$.

\subsection{Continuous time birth-and-death chains and related processes.}

 In the case of a birth-and-death process, the spectral representation (\ref{generalLaplace}) of ${\L}_{f_{0}}(s_1,x)$ can
 be expressed via orthogonal polynomials.
 
  \noindent Let $\pi_0=1$ and $\pi_k={\lambda_0 \lambda_1 \dots \lambda_{k-1} \over \mu_1 \mu_2 \dots \mu_k}$ for all $k \geq 1$.
  Observe that $\pi=[\pi(0), \pi(1), \pi(2), \dots]$ satisfies the detailed balance (reversibility) condition for the process:
  $$\pi_k \lambda_k = \pi_{k+1} \mu_{k+1}~~\text{ for all }~k\geq 0.$$
 Let $P_0(s)\equiv1$, $P_1(s)$, $P_2(s), \dots$
 (where each $P_k(s)$ is a polynomial of $k$th degree)  be constructed recursively as the coordinates of  an eigenvector
 $$P[s]=\left[\begin{array}{c}P_0(s) \\ P_1(s) \\ P_2(s) \\ \vdots \end{array}\right]~~ \text{ satisfying }~~(Q-sI)P[s]=0~.$$
We recall the results of \cite{km1} and \cite{km2}, where it was shown (extending a theorem of J.Favard) that there
 is a probability measure $\mu$ on $(-\infty,0]$ with infinite support such that the polynomials $\{P_k(s)\}_{k=0,1,\dots}$ are orthogonal
 w.r.t. measure $\mu$,
 \begin{eqnarray} \label{orthogonal}
 \int_{(-\infty,0]} P_k(s) P_m(s) d\mu(s)={\delta_{k,m} \over \pi_k}~.
 \end{eqnarray}

 The expression (\ref{generalLaplace}) for ${\L}_{f_0}(s_1,x)$ can be rewritten via the Cauchy transforms w.r.t. the spectral probability measure $\mu$ as follows.
 $$\left((Q-s_1I)^{-1} e_k,e_0 \right)= \int_{(-\infty,0]} {P_k(x) \over  x-s_1} d\mu(x)=C(P_k d\mu)(s_1),$$
 where $C(g d\mu)(s)=\int_{(-\infty,0]} {g(x) \over  x-s} d\mu(x)$ denotes the Cauchy transform (w.r.t. $d\mu$) of $g$.
 Here the function $h(s)$ of (\ref{generalLaplace}) can be expressed as
 $$h(s)= -\int_{(-\infty,0]} {d\mu(x) \over  x-s}~.$$

\subsubsection{Birth-and-death process with equal rates}
Here, we will compute the occupation time density for a birth-and-death process with forward rates $\lambda=1$ and 
reverse rates $\mu=1$, i.e. the process whose generator is a simple Jacobi operator 
$$Q=\left(\begin{array}{ccccc}-r & r & 0 & 0 & \dots \\1 & -2 & 1 & 0 & \dots \\0 & 1 & -2 & 1 & \dots \\0 & 0 & 1 & -2 & \ddots \\\vdots & \vdots & \vdots & \ddots & \ddots\end{array}\right)~, ~~~\text{ where }~r>0.$$
\vskip 0.2 in
{\theorem\label{birth} The zero-state occupation time density for a birth-and-death process with equal rates can be expressed via
modified Bessel functions as follows:
 $$f_0(t,x)=e^{-rt} \delta_0(t-x)+re^{(2-r)x-2t}I_0\left(2\sqrt{(t-x)(t+(r-1)x)}\right)\cdot {\bf 1}_{\{x \leq t\}}$$
 $$+{rt \over \sqrt{(t-x)(t+(r-1)x)}}e^{(2-r)x-2t}I_1\left(2\sqrt{(t-x)(t+(r-1)x)}\right)\cdot {\bf 1}_{\{x \leq t\}}~.$$
 }
 
\begin{proof}
Equation (\ref{generalFourierLaplace}) translates as
\begin{eqnarray*}
{\L}_{\hat{f}_{0}}(s_1,s_2) & = & {1 \over r+s_1-is_2}+{r \over r+s_1-is_2}{\L}_{\hat{f}_{1}}(s_1,s_2)\\
{\L}_{\hat{f}_{1}}(s_1,s_2) & = & {1 \over 2+s_1}+{1 \over 2+s_1}{\L}_{\hat{f}_{0}}(s_1,s_2)+{1 \over 2+s_1}{\L}_{\hat{f}_{2}}(s_1,s_2)\\
{\L}_{\hat{f}_{k}}(s_1,s_2) & = & {1 \over 2+s_1}+{1 \over 2+s_1}{\L}_{\hat{f}_{k-1}}(s_1,s_2)+{1 \over 2+s_1}{\L}_{\hat{f}_{k+1}}(s_1,s_2)~~~~(k=1,2,\dots)
\end{eqnarray*}
where ${\L}_{\hat{f}_{k}}(s_{1},s_{2})$ again denotes the Laplace transform of Fourier transform of $f_k$.
In this recurrence relation, let 
$$l_k(s_1,s_2)={\L}_{\hat{f}_{k}}(s_1,s_2)-{1\over s_1}~.$$
Then $l_k$ satisfy the following recurrence relation, 
\begin{eqnarray*}
l_{k}(s_1,s_2) & = & {1 \over 2+s_1}l_{k-1}(s_1,s_2)+{1 \over 2+s_1}l_{k+1}(s_1,s_2)~~~~(k=1,2,\dots),
\end{eqnarray*}
Solving the characteristic equation,
$$x^2-(2+s_1)x+1=0$$
get 
$$l_k(s_1,s_2)=\alpha_1(s_1,s_2)\left({2+s_1+\sqrt{s_1^2+4s_1} \over 2}\right)^k +\alpha_2(s_1,s_2) \left({2+s_1-\sqrt{s_1^2+4s_1} \over 2}\right)^k~.$$
Observe that ${\L}_{\hat{f}_{k}}(s_1,0)=\int_{[0,+\infty)} \int_{\mathbb{R}}e^{-s_1t} f_k(t,x) dx dt={1\over s_1}$ and
$${\L}_{\hat{f}_{k}}(s_1,s_2)=\int_{[0,+\infty)} \int_{\mathbb{R}}e^{-s_1t+is_2x} f_k(t,x) dx dt \rightarrow 
\int_{[0,+\infty)} \int_{\mathbb{R}}e^{-s_1t+is_2x} \delta_0(x) dx dt ={1\over s_1}~~\text{ as } k \rightarrow \infty~.$$
That is
$$l_k(s_1,s_2) \rightarrow 0 ~~\text{ as } k \rightarrow \infty~.$$
Hence, since $s_1>0$,
\begin{eqnarray*}
l_k(s_1,s_2)=l_0(s_1,s_2) \left({2+s_1-\sqrt{s_1^2+4s_1} \over 2}\right)^k~.
\end{eqnarray*}
Now the first recurrence relation for $\{L_n\}$ can be rewritten as
$${\L}_{\hat{f}_{0}}(s_1,s_2)  =  {1 \over r+s_1-is_2}+{r \over r+s_1-is_2}\left[\Big({\L}_{\hat{f}_{0}}(s_1,s_2)-{1\over s_1}\Big){2+s_1-\sqrt{s_1^2+4s_1} \over 2}+{1\over s_1} \right]$$
Therefore
$${\L}_{\hat{f}_{0}}(s_1,s_2)  = {i \over 2s_1} \cdot {(2-r)s_1+r\sqrt{s_1^2+4s_1} \over s_2+\frac{i}{2}((2-r)s_1+r\sqrt{s_1^2+4s_1})}~.$$
 Once again, using complex integration, we arrive to
$${\L}_{f_{0}}(s_1,x)  = {(2-r)s_1+r\sqrt{s_1^2+4s_1}  \over 2s_1}\exp\left\{-\frac{x}{2}\left((2-r)s_1+r\sqrt{s_1^2+4s_1}\right)\right\}~.$$ 

\vskip 0.2in
 \noindent Observe that one can use the same characteristic equation in order to find the expression for the corresponding
 orthogonal polynomials:
 $$P_k(s)=\left({1 \over 2}-{1 \over \sqrt{s^2+4s}}\right)\left({2+s+\sqrt{s^2+4s} \over 2}\right)^k+\left({1 \over 2}+{1 \over \sqrt{s^2+4s}}\right)\left({2+s-\sqrt{s^2+4s} \over 2}\right)^k~,$$
 where  the spectral measure will satisfy
 $$\int_{(-\infty,0]} {d\mu(x) \over  x-s}={-2 \over (2-r)s+r\sqrt{s^2+4s}}~.$$
 \vskip 0.2in 
 We will now invert the Laplace transform by decomposing ${\L}_{f_{0}}(s_1,x)$ as follows
 $${\L}_{f_{0}}(s_1,x)  = {2-r \over 2}\mathcal{P}_{I}+{r \over 2}\mathcal{P}_{II}+2r\mathcal{P}_{III}~,$$
 where \qquad
 $$\mathcal{P}_{I}=\exp\left\{-\frac{x}{2}(2-r)s_1\right\} \cdot \exp\left\{-\frac{x}{2}r\sqrt{s_1^2+4s_1} \right\},$$
 $$\mathcal{P}_{II}= {s_1 \over \sqrt{s_1^2+4s_1}}\exp\left\{-\frac{x}{2}(2-r)s_1\right\} \cdot \exp\left\{-\frac{x}{2}r\sqrt{s_1^2+4s_1} \right\}$$
  and $$\mathcal{P}_{III}= {1 \over \sqrt{s_1^2+4s_1}}\exp\left\{-\frac{x}{2}(2-r)s_1\right\} \cdot \exp\left\{-\frac{x}{2}r\sqrt{s_1^2+4s_1} \right\}.$$  
  \noindent We will quote a Laplace transform formula (29.3.91) in \cite{as}:
  $$\int_{k}^{\infty} e^{-st} e^{-{1 \over 2}at}I_0\left({1 \over 2}a\sqrt{t^2-k^2} \right) dt 
  = {e^{-k\sqrt{s(s+a)}} \over \sqrt{s(s+a)}},~~~(k \geq 0)~.$$
 First we will find the inverse-Laplace transform of  $\mathcal{P}_{III}$. Taking $s=s_1$, $a=4$ and $k={rx \over 2}$
  in (29.3.91) of \cite{as}, we get
  $$\int_{rx \over 2}^{\infty} e^{-s_1t} e^{-2t}I_0\left(2\sqrt{t^2-\left({rx \over 2}\right)^2} \right) dt 
  = {e^{-{rx \over 2}\sqrt{s_1(s_1+4)}} \over \sqrt{s_1(s_1+4)}}$$  
 Multiplying both sides of the above equation by $\exp\left\{-\frac{(2-r)x}{2}s_1\right\}$, and changing the variable to
 $t:=t+{(2-r)x \over 2}$, obtain
  \begin{eqnarray} \label{p3}
  \mathcal{P}_{III}= \int_{x}^{\infty} e^{-s_1t} e^{(2-r)x-2t}I_0\left(2\sqrt{(t-x)(t+(r-1)x)}\right) dt~.
  \end{eqnarray}
 Therefore, the inverse of $\mathcal{P}_{III}$ is
 $$\mathcal{L}^{-1}(\mathcal{P}_{III})=e^{(2-r)x-2t}I_0\left(2\sqrt{(t-x)(t+(r-1)x)}\right)\cdot {\bf 1}_{\{x \leq t\}}.$$
 We  differentiate $ {\partial \over \partial t}$ and integrate by parts in (\ref{p3}):
 $$\mathcal{L}\left( {\partial \over \partial t} \left[e^{(2-r)x-2t}I_0\left(2\sqrt{(t-x)(t+(r-1)x)}\right) \cdot {\bf 1}_{\{x \leq t\}}\right] \right)
 =\mathcal{P}_{II}-e^{-rx} e^{-s_1x}~.$$
 Hence
 $$\mathcal{L}^{-1}(\mathcal{P}_{II})=
      e^{-rt} \delta_0(t-x) -2e^{(2-r)x-2t}I_0\left(2\sqrt{(t-x)(t+(r-1)x)}\right)\cdot {\bf 1}_{\{x \leq t\}}$$
      $$+{2t-(2-r)x \over \sqrt{(t-x)(t+(r-1)x)}}e^{(2-r)x-2t}I_1\left(2\sqrt{(t-x)(t+(r-1)x)}\right)\cdot {\bf 1}_{\{x \leq t\}}.$$
 In order for us to invert $\mathcal{P}_{I}$, we will need (29.3.96) of \cite{as}, that states the following
  $$\int_{k}^{\infty} e^{-st} {ak \over \sqrt{t^2-k^2}}I_1\left(a\sqrt{t^2-k^2} \right) dt 
  = e^{-k\sqrt{s^2-a^2}} -e^{-ks},~~~(k >0)~.$$ 
 Here we let $s=s_1+2$, $a=2$ and $k={rx \over 2}$, thus obtaining
  $$\mathcal{P}_{I}=e^{-rx} e^{-s_1x}
  +\int_{rx \over 2}^{\infty} e^{-s_1(t+{(2-r)x \over 2})} {rx e^{-2t} \over \sqrt{t^2-\left({rx \over 2}\right)^2}}I_1\left(2\sqrt{t^2-\left({rx \over 2}\right)^2} \right) dt~.$$  
 Once again changing the variable to $t:=t+{(2-r)x \over 2}$, get
  $$\mathcal{P}_{I}=e^{-rx} e^{-s_1x}
  +\int_{x}^{\infty} e^{-s_1t} {rx e^{(2-r)x-2t} \over \sqrt{(t-x)(t+(r-1)x)}}I_1\left(2\sqrt{(t-x)(t+(r-1)x)} \right) dt$$ 
  and
 $$\mathcal{L}^{-1}(\mathcal{P}_{I})=
      e^{-rt} \delta_0(t-x)+{rx e^{(2-r)x-2t} \over \sqrt{(t-x)(t+(r-1)x)}}I_1\left(2\sqrt{(t-x)(t+(r-1)x)} \right)\cdot {\bf 1}_{\{x \leq t\}}~.$$
 We add up all three terms together, thus proving the theorem.
 \end{proof}

\subsubsection{How the occupation time densities are expressed via modified Bessel functions, $I_n(\cdot)$, and 
  the moments of the spectral measure}
   
   Let $m_0, m_1, \dots$ denote the moments of the spectral measure $\mu$, i.e.
   $$m_j=\int_{(-\infty,0]} (-x)^j d\mu(x)~.$$ 
   Consider a case where the spectral measure $\mu$ has bounded support, say $supp(\mu) \subset [-K,0]$.
    Then, for any $z \in \mathbb{C} \setminus (-\infty,0]$ such that $|z|>2K$,
    $${1 \over h(z)}
    ={z \over 1-m_1z^{-1}+m_2z^{-2}-\dots}=z\left(1+\sum_{k=1}^{\infty} (m_1z^{-1}-m_2z^{-2}+\dots)^k\right)$$
    $$=z+m_1-(m_2-m_1^2)z^{-1}+\phi(z^{-1})z^{-2}~.$$
   Recall (\ref{generalLaplace}). Now, we will consider the inverse Laplace transform of a function $F_0(t,x)$, whose
   Laplace transform  ${\L}_{F_0}(s_1,x)={1 \over s_1}\exp\left\{-{x  \over h(s_1)}\right\}$.
   Here for an $a>0$ and $0\leq x \leq t$,
   $$F_0(t,x)={1 \over 2\pi i}\int_{a-i \infty}^{a+i\infty} {1 \over z}\exp\left\{zt-{x \over h(z)} \right\} dz$$
   $$={1 \over 2\pi i}\int_{a-i \infty}^{a+i\infty} {1 \over z}\exp\left\{z(t-x)-m_1x+(m_2-m_1^2)xz^{-1}-x\phi(z^{-1})z^{-2} \right\} dz~.$$
   Here $e^{-x\phi(z^{-1})z^{-2}}=1+\sum_{k=2}^{\infty} v_k(x)z^{-k}$. Thus
   \begin{eqnarray*}
    F_0(t,x) & = & e^{-m_1 x} I_0\left(2\sqrt{(m_2-m_1^2)(t-x)x}\right)\\
    & + & e^{-m_1 x} \sum_{k=2}^{\infty} v_k(x) \left( {t-x \over (m_2-m_1^2)x}\right)^{k \over 2}  I_k\left(2\sqrt{(m_2-m_1^2)(t-x)x}\right)
    \end{eqnarray*}
     by (29.3.81) in \cite{as}.
   
   Another derivation of the same formula may allow further simplification.
   Since $Q$ is the generator for a reversible Markov process, its spectrum is entirely contained inside  $(-\infty,0]$.
   Here the spectrum is a subset of $[-K,0]$. 
   Let $m(z,x)={1 \over 2\pi iz}e^{-x\phi(z^{-1})z^{-2}}$ and for $s \in (-\infty,0]$, $v(s,x)=m_-(s,x)-m_+(s,x)$, where
    $m_-(s,x)=\lim_{\varepsilon \downarrow 0} m(s-i\varepsilon,x)$ and $m_+(s,x)=\lim_{\varepsilon \downarrow 0} m(s+i\varepsilon,x)$. 
   Observe that $m(z,x) \rightarrow 0$ as $|z| \rightarrow \infty$. Solving the Riemann-Hilbert problem
   via Plemelj formula, obtain
   $$m(z,x)={1 \over 2\pi i} \int_{(-\infty,0]} {v(s,x) \over z-s} ds~~\text{ for }~~z\in \mathbb{C}\setminus [-K,0].$$
   Now,
   $$F_0(t,x)=\int_{a-i \infty}^{a+i\infty} m(z,x) \exp\left\{z(t-x)-m_1x+(m_2-m_1^2)xz^{-1}\right\} dz$$
   $$= \int_{(-\infty,0] } \left[ {1 \over 2\pi i} \int_{a-i \infty}^{a+i\infty} {1 \over z-s} \exp\left\{z(t-x)-m_1x+(m_2-m_1^2)xz^{-1}\right\} dz \right] v(s,x) ds$$ 
   $$=  e^{-m_1 x}\int_{(-\infty,0] } \sum_{k=0}^{\infty} s^k \left[ {1 \over 2\pi i} \int_{a-i \infty}^{a+i\infty} {1 \over z^{k+1}} \exp\left\{z(t-x)+(m_2-m_1^2)xz^{-1}\right\} dz \right] v(s,x) ds$$
     $$=  e^{-m_1 x} \sum_{k=0}^{\infty} \int_{(-\infty,0] }s^k v(s,x) ds \left( {t-x \over (m_2-m_1^2)x}\right)^{k \over 2}
     I_k\left(2\sqrt{(m_2-m_1^2)(t-x)x}\right)$$
     by (29.3.81) in \cite{as}.
   
   \noindent
   Observe that since $\phi(z) z^{-2}$ has no poles at $\infty$, for $k \geq 1$,
    $$\int_{(-\infty,0] } s^k v(s,x) ds=\oint_{\gamma_+ \cup \gamma_-} z^{k-1}e^{-x \phi(z^{-1}) z^{-2}} dz~,$$
    where $\gamma_+$ connects the origin to $\infty$ right above the negative half-line, and 
    $\gamma_-$ connects $\infty$ to the origin barely bellow the negative half-line. Observe that for $k=1$ the above integral is zero.
   
   Here $f_0(t,x)=-{\partial \over \partial x} F_0(t,x)$.

\section{Conclusions}

 Observe that in general, an argument in Deift \cite{deift} (Section 2) shows that if a process is time reversible (i.e. satisfies a detailed balance condition) 
 with a bounded generator then there exist a unique (spectral) probability measure $\mu$ with compact support $supp(\mu) \subset \mathbb{R}$ such that
 $$h(s_1)=-\left(e_0,(Q-s_1I)^{-1}e_0\right)=-\int_{\mathbb{R}} {d\mu(x) \over x-s_1} $$  
 As a conclusion, let us list some of the open problems and directions for further research the authors are working on.
\begin{itemize}
 \item What properties of $d\mu$ would allow for the inversion of the Laplace transform ${\L}_{f_0}(s_1,x)$
  to be expressed explicitly via modified Bessel functions $I_n$?
 \item Exploring occupation time densities for a wider class of time reversible stochastic processes.
 \item Interpreting reinforced processes studied in  Kovchegov \cite{kov} as occupation time driven processes.
\end{itemize}

\section*{Acknowledgment}
The authors wish to thank R.Burton and M.Ossiander for sharing thoughts on the subject of this paper.

\bibliographystyle{amsplain}

\end{document}